\documentclass{article}
\usepackage[utf8]{inputenc}
\usepackage{geometry}                		
\usepackage{graphicx}						
\usepackage{mathtools}
\usepackage{mathrsfs}
\usepackage[all]{xy}
\usepackage{amsmath}
\usepackage{MnSymbol}

\usepackage{wasysym}
\usepackage{tikz,tikz-cd}
\newtheorem{theorem}{Theorem}

\title{A new proof of Huber's theorem on differential geometry in the large}
\author{Chen Zhou}
\date{ }

\begin{document}

\maketitle

\begin{abstract}
In this paper we give a new, and shorter, proof of Huber's theorem (Theorem 13 in \cite{Huber1958OnSF}) which affirms that for a connected open Riemann surface endowed with a complete conformal Riemannian metric, if the negative part of its Gaussian curvature has finite mass, then the Riemann surface is homeomorphic to the interior of a compact surface with boundary, and thus it has finite topological type. We will also show that such Riemann surface is parabolic. 
	
\end{abstract}

\section{Introduction}
The connection between geometry and topology has been an attractive topic for a long time. We have so many interesting results about it now, such as the classical Gauss-Bonnet formula and its extensions to surfaces with boundaries and to higher dimensions. 

The developments of comparison theorems in Riemannian geometry enrich the study a lot in the case when the manifold has non-positive or non-negative curvature, and in such cases we do not need to assume the compactness of the manifold, but rather that it carries a complete metric. For example, when the manifold has a complete Riemannian metric of class $\mathcal{C}^\infty$ with sectional curvature bounded below by zero, we have the famous ``soul theorem"\cite{10.2307/1970819}. A theorem of Gromov (1.5 in \cite{10.4310/jdg/1214434488}) even gives an estimate of the number of generators for the fundamental group of such manifold. If we realise the condition to case when the curvature of the manifold is ``not very negative at infinity", Abresch showed in  \cite{Abresch1985} and \cite{ASENS_1987_4_20_3_475_0} that the Betti numbers of the manifold are finite, by generalizing the idea of Gromov in \cite{Gromov1981}. In fact, he showed in \cite{ASENS_1987_4_20_3_475_0} that all the critical points of the distance function on such kind of manifold are contained in a compact subset. Then by the Isotopy Lemma, the manifold should be homeomorphic to the interior of a compact manifold with boundary.

In the two (real) dimensional case, Huber proved a very interesting result (Theorem 13 in \cite{Huber1958OnSF}). His theorem does not require a restriction on the bound of the curvature, but just assumes the negative part of the Gaussian curvature has finite mass. In order to describe his theorem, we first need to fix some notation. Let $M$ be a Riemann surface and $ds^2$ be a conformal Riemannian metric of class $\mathcal{C}^{\infty}$ on $M$. This metric has an expression
$$ds^2=e^{2u(z)}|dz|^2$$
in local coordinates. Its Gaussian curvature  can be expressed as 
$$K=-e^{-2u}\Delta u,$$
where 
$$\Delta=4\,\frac{\partial^2}{\partial z\partial\bar{z}}$$ 
is the (local) Laplace operator. Denote by $d\sigma$ the area element on $M$. It can be written as 
$$d\sigma=\frac{i}{2}\,e^{2u}\,dz\wedge d\bar{z}$$
in the above coordinates. Therefore locally we have 
$$K d\sigma=-\frac{i}{2}\,\Delta u\,dz\wedge d\bar{z}.$$
If we define 
$$K^-:=\max\{-K, 0\}$$
and 
$$\Delta^+u:=\max\{\Delta u, 0\} \;\;\;\;\;\;\;\; \Delta^-u:=\max\{-\Delta u, 0\},$$
then similarly we have
$$K^-=e^{-2u}\Delta^+ u$$
and
$$K^- d\sigma=\frac{i}{2}\,\Delta^+u\,dz\wedge d\bar{z}$$
locally. Huber's theorem then can be stated as:
\begin{theorem}[Huber \cite{Huber1958OnSF}]\label{th1}
Let $M$ be a connected open Riemann surface and let $ds^2$ be a complete conformal Riemannian metric of class $\mathcal{C}^{\infty}$ on $M$. 
If the integral 
$$\int_MK^-d\sigma$$
is finite, then  $M$ is homeomorphic to the interior of a compact surface with boundary.
\end{theorem}
The original proof of the above theorem given by Huber contains many long computations in complex analysis and it takes great effort for people to go through it. However, the developments of geometry and topology \cite{Abresch1985} and \cite{ASENS_1987_4_20_3_475_0} emerging after the publication of Huber's paper provide many new ideas which enable us to give a shorter, and more conceptual, proof of his theorem. But in Abresch's papers \cite{Abresch1985} and \cite{ASENS_1987_4_20_3_475_0}, he assumes that the curvature is ``not very negative at infinity", so we should first reduce our problem in Theorem \ref{th1} to that case. This can be done with the help of the potential theory on Riemann surfaces (see section 7-13 of Chapter V in \cite{nla.cat-vn1098973} and see also Appendix I of \cite{Sario2013value}). Nevertheless, we may lose the smoothness of the Riemannian metric during this process; the Riemannian metric of class $\mathcal{C}^{\infty}$ may then become a Riemannian metric of class $\mathcal{C}^2$. But noting that all the proofs in \cite{Abresch1985} and \cite{ASENS_1987_4_20_3_475_0} are also valid for a smooth manifold endowed with a Riemannian metric of class $\mathcal{C}^2$, a new proof of Huber's theorem can be given.


\section{A new proof of Theorem \ref{th1}}

Before giving the proof, let us first recall some concepts about the potential theory on Riemann surfaces (see Appendix I of \cite{Sario2013value} for more details, see also section 1 of \cite{Napier1995StructureTF}). An open Riemann surface (or more generally, an open Riemannian manifold) $M$ is called \textit{hyperbolic} if it admits a positive Green's function $G$. This is equivalent to the existence of a non-constant positive continuous superharmonic function on $M$. A Riemann surface which is not hyperbolic is called \textit{parabolic}.

\textit{Proof of Theoren \ref{th1}}. For the Riemann surface $M$ in Theoerm \ref{th1}, let $D_2$ be a region on $M$ that is biholomorphic to $\{z: |z|\leqslant 2\}$ in the complex plane and $D_1$ be the subset of $D_2$ that corresponds to $\{z: |z|\leqslant 1\}$. Consider the Riemann surface 
$$M'=M-D_1.$$ 
It is a hyperbolic Riemann surface. In fact, letting $S_1$ be the subset of $D_2$ that corresponds to $\{z: |z|=1\}$ and $S_2$ the subset that corresponds to $\{z: |z|=2\}$, we know there exists a real-valued continuous function $h$ that equals $1$ on $S_1$, equals $0$ on $S_2$ and is harmonic on the region $E$ between $S_1$ and $S_2$. By the third characterization of hyperbolicity on page 813 in \cite{Napier1995StructureTF}, $M'$ is a hyperbolic Riemann surface. 

Let $G'$ be the Green's function on $M'$. We denote the restriction of $ds^2$ on $M'$ still by $ds^2$ and we do the same thing for $d\sigma$ and $K^-$ as well. Then $\int_MK^-d\sigma$ being finite implies the following integral converges:
$$f(p):=\int_{M'} G'(p, q)\,K^-(q)\,d\sigma(q).$$
In fact, for every fixed point $p_0$ in $M'$, we can find a coordinate neighborhood $U_{p_0}$ of $p_0$ such that $G(p_0, q)$ is equal to a constant $C$ times $(h(z)-\log |z_0-z|)$. ( Here $z_0$ is the coordinate of $p_0$ and $z$ is the coordinate of $q$ while $h(z)$ is a continuous function on $U_{p_0}$.) Thus $G'(p_0, q)\,K^-(q)$ is an integrable function of $q$ in $U_{p_0}$. On the other hand,  since $G'(p_0, q)$ is a bounded function of $q$ outside $U_{p_0}$, the above integral also converges outside $U_{p_0}$.

We can show the function $f(\cdot)$ defined above is a $\mathcal{C}^2$ function. In fact, since the original metric $ds^2$ is smooth, its Gaussian curvature $K$ is a smooth function. The function $K^-$ is defined as 
$$K^- :=\max\{-K, 0\},$$ so it is a locally Holder continuous function for some exponent $\alpha$. 
By the elliptic regularity, the function
$$f(p):=\int_{M'}G'(p, q)K^-(q)d\sigma(q)$$
is a locally $\mathcal{C}^{2,\alpha}$ function on $M'$.  (See Lemma 4.2 in \cite{MR737190}, and note that the Green's function locally equals a harmonic function plus the fundamental solution of Laplace's equation.)

We can extend $f$ to a $\mathcal{C}^{2}$ function on $M$, thus the metric $ds_0^2=e^{2f}ds^2$ defined on $M$ is a conformal metric of class $\mathcal{C}^{2}$, and locally can be written as
$$ds_0^2=e^{2f}ds^2=e^{2V(z)}|dz|^2,$$
where 
$$V(z)=f(z)+u(z).$$
If we denote the (global) Laplace operator associated with the Riemannian metric $ds^2$ by $\Delta'$, then on $M'$ we should have 
$$\Delta' f=-K^-.$$
Noting that $\Delta'$ can be locally expressed as 
$$\Delta'=4\,e^{-2u}\,\frac{\partial^2}{\partial z\partial\bar{z}}=e^{-2u}\Delta,$$
we get
$$\Delta f=-e^{2u}K^-=-\Delta^+u$$
on $M'$, and hence on $M'$, the Gaussian curvature of the metric $ds_0^2$ can be locally expressed as 
$$-e^{-2V}\Delta V=-e^{-2V}(\Delta u-\Delta^+u)=e^{-2V}\Delta^-u,$$
 which is non-negative. Since $D_1$ is compact, we have the Gaussian curvature of $ds_0^2$ is non-negative on $M$ outside a compact set. 

Since we assume that $(M, ds^2)$ is complete, thus for any smooth curve $\gamma$ going to infinity, its length in the metric $ds^2$ will be infinite. Noting that $ds_0^2=e^{2f}ds^2$ and $f$ is non-negative outside a compact subset of $M$, the length of the above curve $\gamma$ in the metric $ds_0^2$  will also be infinite. This implies $(M,ds_0^2)$ is also complete. 
 Since $(M, ds_0^2)$ is complete and the Gaussian curvature of the metric $ds_0^2$ on $M$ is non-negative outside a compact subset, we can then use the results of Abresch in \cite{Abresch1985} 
 and \cite{ASENS_1987_4_20_3_475_0} (although Abresch assumes that the manifold is a smooth Riemannian manifold in his papers, his proof is valid for any smooth manifold endowed with a Riemannian metric of class $\mathcal{C}^2$): 

   Let $d_0(\cdot, \cdot)$ be the distance function associated to the metric $ds_0^2$. Fix a point $p$ on $M$ and define 
 $$d_p(\cdot):=d_0(p, \cdot).$$
 A point $q$ is called a \textit{critical point} of $d_p(\cdot)$ if for any tangent vector $v \in T_q M$, there is a minimal geodesic segment $\gamma$ which joins $\gamma(0)=q$ to $p$ such that the angle between $\gamma'(0)$ and $v$ is less than or equal to $\frac{\pi}{2}$. We can prove that all the critical points of the function $d_p$ are contained in a compact subset of $M$ as follows. Suppose this is not true. Then for any fixed $L>1$ and any natural number $k$, we can find $k$ critical points of $d_p(\cdot)$, which we denote by $\{p_1,\cdots, p_k\}$,  such that 
 $$d_p(p_{j-1})\geqslant L\,d_p(p_j).$$
 But by 2.3 Lemma (i) on page 481 in \cite{ASENS_1987_4_20_3_475_0}, $k$ is bounded above by a constant, so we get a contradiction.
 
 Now let 
$$B_0(p; r):=\{q\in M: d_p(q)=d_0(p, q)< r\}.$$ We have the following Isotopy Lemma (for the proof, see Lemma 1.1 in \cite{Gromov1981}, and note that the Iostopy Lemma is even true for Alexandrov spaces as pointed in Lemma 1 in \cite{10.2307/2160867}):

\textbf{Isotopy Lemma}
\textit{Suppose that the closure of ($B_0(p; r_2)-B_0(p; r_1))$ does not contain any critical points of $d_p(\cdot)$. Then there is an isotopy of $M$ which sends $B_0(p; r_2)$ to $B_0(p; r_1)$ and fixes the points outside a neighborhood of $B_0(p; r_2)$.}

By the above Isotopy Lemma, $M$ is homeomorphic to the interior of $P$, where $P (\subset M)$ is a compact bordered Riemann surface. This proves the theorem.\\
 
\section{Parabolicity of $M$}
 We can also prove the following (Theorem 15 in \cite{Huber1958OnSF}):
\begin{theorem}[Huber \cite{Huber1958OnSF}]\label{th2}
Let $M$ be a connected open Riemann surface and let $ds^2$ be a complete conformal Riemannian metric of class $\mathcal{C}^{\infty}$ on $M$. 
If the integral 
$$\int_MK^-d\sigma$$
is finite, then $M$ parabolic.
\end{theorem}

\textit{Proof of Theorem \ref{th2}}. If not, suppose that $M$ is hyperbolic. Define
$$f_1(p):=\int_{M}G(p, q)K^-(q)d\sigma(q).$$
As above, we know $f_1$ is a $\mathcal{C}^2$ function on $M$. Thus $ds_1^2=e^{2f_1}ds^2$ defines a $\mathcal{C}^2$ conformal metric on $M$ with non-negative Gaussian curvature. Fix a point $p$ on $M$. Let $d_1(\cdot, \cdot)$ be the distance function associated to the metric $ds_1^2$. Define 
$$B_1(p; r):=\{q\in M: d_1(p, q)< r\}.$$ 
Let $A_1$ be the area function associated to the metric $ds_1^2$ and denote 
$$A_1(p; r):=A_1(B_1(p; r)).$$
By the Bishop-Gromov theorem (see 1.9.4 in \cite{karcher1987riemannian} for Riemannian manifolds and Theorem 10.6.6 in \cite{burago2001course} for Alexandrov spaces), we have 
$$A_1(p; r)\leqslant 2\pi r^2$$
as $r$ goes to infinity. Then by the following result of Cheng and Yau \cite{Cheng1975DifferentialEO}, $M$ is parabolic. This contradicts the assumption that $M$ is hyperbolic.

\textbf{Theorem (Cheng and Yau \cite{Cheng1975DifferentialEO})}
\textit{Let $M$ be a complete Riemannian manifold and let $p$ be a point on $M$. If the volume $V_r$ of a geodesic ball centered at $p$ with radius $r$ satisfies the relation 
$$\liminf_{r\to \infty} \frac{V_r}{r^2}<\infty,$$
then $M$ parabolic.}

For the proof of the above theorem, see \cite{Cheng1975DifferentialEO}. One can also find another proof of this theorem in \cite{Grigor'yan_1987}, see Corollary 3 in \cite{Grigor'yan_1987}.

\section*{Acknowledgement}
I would like to thank Professor Mohan Ramachandran for giving me this interesting problem and also providing many enlightening and helpful suggestions about it. His instruction through both face-to-face conversations and emails is crucial for me to find a new proof of Huber's theorem.

{
 \bibliographystyle{ieee}
\bibliography{Cmm}

\begin{thebibliography}{10}\itemsep=-1pt

\bibitem{Abresch1985}
U.~Abresch.
\newblock Lower curvature bounds, toponogov's theorem, and bounded topology.
\newblock {\em Annales scientifiques de l'École Normale Supérieure},
  18(4):651--670, 1985.

\bibitem{ASENS_1987_4_20_3_475_0}
U.~Abresch.
\newblock Lower curvature bounds, {Toponogov's} theorem, and bounded topology.
  {II}.
\newblock {\em Annales scientifiques de l'\'Ecole Normale Sup\'erieure}, Ser.
  4, 20(3):475--502, 1987.

\bibitem{burago2001course}
D.~Burago, Y.~Burago, and S.~Ivanov.
\newblock {\em A Course in Metric Geometry}.
\newblock Graduate Studies in Mathematics. American Mathematical Society, 2001.

\bibitem{10.2307/1970819}
J.~Cheeger and D.~Gromoll.
\newblock On the structure of complete manifolds of nonnegative curvature.
\newblock {\em Annals of Mathematics}, 96(3):413--443, 1972.

\bibitem{Cheng1975DifferentialEO}
S.-Y. Cheng and S.-T. Yau.
\newblock Differential equations on riemannian manifolds and their geometric
  applications.
\newblock {\em Communications on Pure and Applied Mathematics}, 28:333--354,
  1975.

\bibitem{MR737190}
D.~Gilbarg and N.~S. Trudinger.
\newblock {\em Elliptic partial differential equations of second order}, volume
  224 of {\em Grundlehren der Mathematischen Wissenschaften [Fundamental
  Principles of Mathematical Sciences]}.
\newblock Springer-Verlag, Berlin, second edition, 1983.

\bibitem{Grigor'yan_1987}
A.~A. Grigor'yan.
\newblock On the existence of positive fundamental solutions of the laplace
  equation on riemannian manifolds.
\newblock {\em Mathematics of the USSR-Sbornik}, 56(2):349, feb 1987.

\bibitem{10.4310/jdg/1214434488}
M.~Gromov.
\newblock {Almost flat manifolds}.
\newblock {\em Journal of Differential Geometry}, 13(2):231 -- 241, 1978.

\bibitem{Gromov1981}
M.~Gromov.
\newblock Curvature, diameter and betti numbers.
\newblock {\em Commentarii mathematici Helvetici}, 56:179--195, 1981.

\bibitem{Huber1958OnSF}
A.~Huber.
\newblock On subharmonic functions and differential geometry in the large.
\newblock {\em Commentarii Mathematici Helvetici}, 32:13--72, 1958.

\bibitem{karcher1987riemannian}
H.~Karcher.
\newblock {\em Riemannian Comparison Constructions}.
\newblock Rheinische Friedrich-Wilhelms-Universit{\"a}t Bonn,
  Sonderforschungsbereich 256, Nichtlineare Partielle Differentialgleichungen.
  SFB 256, 1987.

\bibitem{10.2307/2160867}
L.-K. Koh.
\newblock Betti numbers of alexandrov spaces.
\newblock {\em Proceedings of the American Mathematical Society},
  122(1):247--252, 1994.

\bibitem{Napier1995StructureTF}
T.~Napier and M.~Ramachandran.
\newblock Structure theorems for complete k{\"a}hler manifolds and applications
  to lefschetz type theorems.
\newblock {\em Geometric \& Functional Analysis GAFA}, 5:809--851, 1995.

\bibitem{nla.cat-vn1098973}
L.~Sario and M.~Nakai.
\newblock {\em Classification theory of Riemann surfaces}.
\newblock Springer-Verlag Berlin, New York, 1970.

\bibitem{Sario2013value}
L.~Sario and K.~Noshiro.
\newblock {\em Value Distribution Theory}.
\newblock The university series in higher mathematics. Springer New York, 2013.

\end{thebibliography}

}

\end{document}